\documentclass[12pt]{article}
\usepackage{amsfonts}
\usepackage{amsmath}
\newtheorem{theorem}{Theorem}

\def\qed{\quad \rule{1.5ex}{1.5ex} \bigskip}
\begin{document}
\title{Multiple orthogonal polynomials associated with Macdonald functions}
\author{Walter Van Assche \thanks{Research Director of the
Belgian National Fund for Scientific Research. Part of this work
is supported by INTAS 93-219ext and FWO research project 
G.0278.97} \\
Semyon B. Yakubovich \thanks{
On leave from the Belarusian State University, Minsk. Supported by junior fellowship F/96/13 of the
Research Counsel of the Katholieke Universiteit Leuven, BELGIUM}\\
 Katholieke Universiteit Leuven}
\date{}
\maketitle

\begin{abstract}
We consider multiple orthogonal polynomials corresponding to two
Macdonald functions (modified Bessel functions of the second kind),
with emphasis on the polynomials on the diagonal of the Hermite-Pad\'e
table. We give some properties of these polynomials: differential
properties, a Rodrigues type formula and  explicit formulas for
the third order linear recurrence relation.
\end{abstract}

\section{Macdonald functions}

In this paper we will investigate certain polynomials satisfying orthogonality properties
with respect to weight functions related to Macdonald functions $K_\nu(z)$. These
Macdonald functions are modified Bessel functions of the second kind and satisfy
the differential equation
$$  z^2u'' + zu' - (z^2+\nu^2)u = 0, $$
for which they are the solution that remains bounded as $z$ tends to infinity on the real line. They
can be given by the following integral representations
\begin{eqnarray}
  K_\nu(z) & = & \left( \frac{\pi}{2z} \right)^{1/2} \frac{e^{-z}}{\Gamma(\nu+1/2)}
      \int_0^\infty e^{-t} t^{\nu-1/2} (1+\frac{t}{2z})^{\nu-1/2}\, dt  \label{eq:poisson} \\
      & = & \frac12 \left( \frac{z}{2} \right)^\nu 
  \int_0^\infty \exp(-t-\frac{z^2}{4t})
     t^{-\nu-1}\, dt  \label{eq:sommerfeld}
\end{eqnarray}
(see, e.g., The Bateman Manuscript Project \cite[Chapter VII]{bateman}, Nikiforov and Uvarov \cite[pp. 223--226]{nikiuv}).
Useful properties are
\begin{eqnarray}
   -2K'_\nu(z) & = & K_{\nu-1}(z)+K_{\nu+1}(z),  \label{eq:diff1} \\
   - \frac{2\nu}{z} K_\nu(z) & = & K_{\nu-1}(z) - K_{\nu+1}(z). \label{eq:diff2}
\end{eqnarray}
These functions have the asymptotic behavior
\begin{equation}  \label{eq:nearinf}
   K_\nu(z) = \left( \frac{\pi}{2z} \right)^{1/2} e^{-z} [1+ O(1/z)], \qquad z \to \infty, 
\end{equation}
and near the origin
\begin{equation}  \label{eq:near0}
 \begin{cases}
  z^\nu K_\nu(z) = 2^{\nu-1} \Gamma(\nu) + o(1) & \text{if $z \to 0$} \\
  K_0(z) = -\log z + O(1),  & \text{if $z \to 0$}
  \end{cases}   
\end{equation}
The moments of these functions are given explictly by
\begin{equation} \label{eq:Kmoments}
  \int_0^\infty x^\mu K_\nu(x)\,dx = 2^{\mu-1} \Gamma(\frac{\mu+\nu+1}2)
     \Gamma(\frac{\mu-\nu+1}2), \qquad \Re(\mu \pm \nu+1) > 0.
\end{equation}

We will use as weight functions the scaled Macdonald functions
\begin{equation}
  \rho_\nu(x) = 2 x^{\nu/2} K_\nu(2\sqrt{x}), \qquad x > 0. 
\end{equation}
which behave like $\exp(-2\sqrt{x})$ as $x \to \infty$.
From the  properties (\ref{eq:diff1}) and (\ref{eq:diff2}) one easily
obtains the differential properties
\begin{eqnarray}
  (x^{-\nu} \rho_\nu)' & = & - x^{-(\nu+1)}\rho_{\nu+1}, \label{eq:rhoup}\\
  \rho_{\nu+1}' & = & - \rho_\nu , \label{eq:rhodown}
\end{eqnarray}
so that $\rho_\nu$ satisfies the second order linear differential
equation $[x^{\nu+1} (x^{-\nu} \rho_\nu)']' = \rho_\nu$. By a simple change of variable in (\ref{eq:Kmoments})
we can find the moments of $\rho_\nu$ which are given by
\begin{equation}  \label{eq:rhomoment}
  \int_0^\infty x^\mu \rho_\nu(x)\, dx = \Gamma(\mu+\nu+1) \Gamma(\mu+1),
\end{equation}
in particular we see that the $n$th moment of $\rho_0$ is given by $(n!)^2$.

In \cite{open} A. P. Prudnikov formulated as an open problem the construction of the orthogonal polynomials associated with the weight
function $\rho_{\nu}$  (and for more general weights, known as ultra-exponential weight functions). In the present paper we will show that it is more natural to investigate multiple orthogonal polynomials for two
Macdonald weights $\rho_\nu$ and $\rho_{\nu+1}$ since for these multiple
orthogonal polynomials one has nice differential properties, a Rodrigues
formula, and an explicit recurrence relation.

\setcounter{equation}{0}
\section{Multiple orthogonal polynomials}
We will investigate two types of multiple orthogonal polynomials for
the system of weights $(\rho_\nu,\rho_{\nu+1})$ $(\nu \geq 0)$ on $[0,\infty)$ with an additional factor $x^\alpha$.
Let $n, m \in \mathbb{N}$ and $\alpha > -1$, then the
{\it type 1 polynomials} $(A_{n,m}^\alpha,B_{n,m}^\alpha)$ are such that 
$A_{n,m}^\alpha$ is of degree $n$ and
$B_{n,m}^\alpha$ is of degree $m$, and they satisfy the orthogonality conditions
\begin{equation}   \label{eq:type1ortho}
 \int_0^\infty [A_{n,m}^\alpha(x)\rho_\nu(x) + 
   B_{n,m}^\alpha(x) \rho_{\nu+1}(x)] x^{k+\alpha}\, dx = 0,
  \qquad k=0,1,2,\ldots,n+m. 
\end{equation}
This gives $n+m+1$ linear homogeneous equations for the $n+m+2$ unknown coefficients of the polynomials $A_{n,m}^\alpha$ and $B_{n,m}^\alpha$,
 so that we
can find the type 1 polynomials up to a multiplicative factor, which we will fix later. It will be convenient to use the notation
\[   q_{n,m}^\alpha(x) = A_{n,m}^\alpha(x)\rho_\nu(x) + 
   B_{n,m}^\alpha(x) \rho_{\nu+1}(x). \]
{\it Type 2 polynomials} $p_{n,m}^\alpha$ are such that $p_{n,m}^\alpha$ 
is of degree
$n+m$ and satisfies the multiple orthogonality conditions
\begin{eqnarray}
  \int_0^\infty p_{n,m}^\alpha(x) \rho_\nu(x) x^{k+\alpha} \, dx = 0 & &
  \quad k=0,1,2,\ldots,n-1, \label{eq:type2ortho1} \\
 \int_0^\infty p_{n,m}^\alpha(x) \rho_{\nu+1}(x) x^{k+\alpha} \, dx = 0 & &
  \quad k=0,1,2,\ldots,m-1. \label{eq:type2ortho2}
\end{eqnarray}
This means that we distribute the $n+m$ orthogonality conditions
over the two weights $x^\alpha \rho_\nu$ and $x^\alpha \rho_{\nu+1}$.
This gives $n+m$ linear homogeneous equations for the $n+m+1$ unknown
coefficients of $p_{n,m}^\alpha$. For type 2 polynomials we will consider
monic polynomials, thereby fixing the leading coefficient to be $1$ and
leaving $n+m$ coefficients to be determined from 
(\ref{eq:type2ortho1})--(\ref{eq:type2ortho2}).

Multiple orthogonal polynomials are related to Hermite-Pad\'e simultaneous
rational approximation of a system of Markov functions near infinity, (see, e.g., Nikishin and Sorokin \cite{nikisor}
and Aptekarev \cite{aptek}). In the present situation the functions
to be approximated are
\[   f_1(z) = \int_0^\infty \frac{x^\alpha \rho_\nu(x)}{z-x} \, dx,
  \quad
    f_2(z) = \int_0^\infty \frac{x^\alpha \rho_{\nu+1}(x)}{z-x} \, dx. \]
Type 1 Hermite-Pad\'e approximation (Latin type) consists of finding
polynomials $A_{n,m}$, $B_{n,m}$, and $C_{n,m}$ such that
\[  A_{n,m}(z) f_1(z) + B_{n,m}(z) f_2(z) - C_{n,m}(z) = O(z^{-n-m-2}) 
   \qquad z \to \infty,  \]
and type 2 Hermite-Pad\'e approximation (German type) is simultaneous rational approximation of $(f_1,f_2)$ with a common denominator
$p_{n,m}$, i.e., one wants polynomials $p_{n,m}$, $R_{n,m}$ and $S_{n,m}$
such that 
\begin{eqnarray*}
  p_{n,m}(z)f_1(z) - R_{n,m}(z) = O(z^{-n-1}), & & \quad z \to \infty, \\
    p_{n,m}(z)f_2(z) - S_{n,m}(z) = O(z^{-m-1}), & & \quad z \to \infty.
\end{eqnarray*}
The polynomials $A_{n,m}, B_{n,m}$ then are precisely the 
type 1 polynomials $A_{n,m}^\alpha$  and $B_{n,m}^\alpha$ and the polynomial $p_{n,m}$ is the type 2 polynomial $p_{n,m}^\alpha$. The numerator 
$C_{n,m}$ is then given by
\[    C_{n,m}(z) = \int_0^\infty \left[ \frac{A_{n,m}^\alpha(z)-  A_{n,m}^\alpha(x)}{z-x} \rho_\nu(x) + \frac{B_{n,m}^\alpha(z)-  B_{n,m}^\alpha(x)}{z-x} \rho_{\nu+1}(x) \right]
   x^\alpha \, dx, \]
and for type 2 approximation the numerators are
\begin{eqnarray*}
   R_{n,m}(z) & = & \int_0^\infty \frac{p_{n,m}^\alpha(z)-         p_{n,m}^\alpha (x)}{z-x}    \rho_\nu(x) x^\alpha\, dx , \\
   S_{n,m}(z) & = & \int_0^\infty \frac{p_{n,m}^\alpha(z)-            p_{n,m}^\alpha (x)}{z-x} \rho_{\nu+1}(x) x^\alpha\, dx . 
\end{eqnarray*}

Observe that we can write the system $(f_1,f_2)$ as
\[  f_1(z) = \int_0^\infty f(x) \frac{\rho_{\nu+1}(x) x^\alpha}{z-x}\, dx , \quad  f_2(z) = \int_0^\infty  \frac{\rho_{\nu+1}(x) x^\alpha}{z-x}\, dx, \]
where $f$ is itself a Markov function 
\[   f(x) = \frac{\rho_{\nu}(x)}{\rho_{\nu+1}(x)} =
     \frac{1}{\pi^2} \int_0^\infty \frac{s^{-1}\,ds}{(x+s) [J_{\nu+1}^2(2\sqrt{s}) + Y_{\nu+1}^2(2\sqrt{s})]},      \]
which follows from a result by Ismail \cite{ismail}. 
Note however that $f$ is a Markov function of a positive measure
for which not all the moments exist. Nevertheless we can still
say that $(f_1,f_2)$ is a {\it Nikishin system}, which guarantees that the polynomials
$A_{n,m}^\alpha$, $B_{n,m}^\alpha$ and $p_{n,m}^\alpha$ all can be
computed and their degrees are exactly $n$, $m$, and $n+m$ respectively.
Furthermore the zeros of $q_{n,n}^\alpha$, $q_{n+1,n}^\alpha$,
 $p_{n,n}^\alpha$ and $p_{n+1,n}^\alpha$ will
all be on $(0,\infty)$. 

\setcounter{equation}{0}
\section{Differential properties}
In this section we will give some differential properties for the type 1
and type 2 multiple orthogonal polynomials. We will only consider the
multiple orthogonal polynomials on the diagonal or close to the diagonal,
i.e., when $n=m$ or $n=m+1$. These are the most natural and they satisfy
interesting relations.

\begin{theorem}  \label{thm:type2}
For the type 2 multiple orthogonal polynomials 
we have for every $\alpha > -1$
\begin{equation}  \label{eq:type2dif}
  \frac{d}{dx} p_{n,n}^\alpha(x) = 2n p_{n,n-1}^{\alpha+1}(x), \quad
  \frac{d}{dx} p_{n,n-1}^\alpha(x) = (2n-1) p_{n-1,n-1}^{\alpha+1}(x).
\end{equation}
\end{theorem}
This means that the differential operator acting on type 2 multiple orthogonal polynomials lowers the degree by one and raises the $\alpha$
by one. This should be compared with the corresponding differential
property
\[   \frac{d}{dx} L_n^\alpha(x) = -L_{n-1}^{\alpha+1}(x)  \]
for Laguerre polynomials (see, e.g., Szeg\H{o} \cite[Eq.~(5.1.14) on p.~102]{szego}). The normalizing constant is different here since Laguerre
polynomials are not monic polynomials.  

{\bf Proof:}
We begin by using (\ref{eq:type2ortho1})
\[   \int_0^\infty p_{n,n}^\alpha(x) \rho_\nu(x) x^{k+\alpha}\, dx = 0,
  \qquad k=0,1,\ldots,n-1, \]
and then use (\ref{eq:rhodown}) to find
\[  -\int_0^\infty p_{n,n}^\alpha(x) \rho_{\nu+1}'(x) x^{k+\alpha}\, dx = 0,
   \qquad k=0,1,\ldots,n-1. \]
Integration by parts then gives
\[  \left. - p_{n,n}^\alpha(x) \rho_{nu+1}(x) x^{k+\alpha} \right|_0^\infty 
  + \int_0^\infty \left[ p_{n,n}^\alpha(x)x^{k+\alpha} \right]' 
    \rho_{\nu+1}(x)\, dx = 0, \qquad k=0,1,\ldots,n-1. \]
The integrated terms vanish for $k\geq 1$ and $\alpha > -1$ because of
the behavior of $\rho_{\nu +1}$ near $0$ and $\infty$, given in
(\ref{eq:nearinf}) and (\ref{eq:near0}). Now work out the differentiation
of the product, then we have
\begin{multline*}
  \int_0^\infty [p_{n,n}^\alpha(x)]' x^{k+\alpha} \rho_{\nu+1}(x)\, dx \\
   +\ (k+\alpha) \int_0^\infty p_{n,n}^\alpha(x) x^{k+\alpha-1} \rho_{\nu+1}(x)\, dx = 0, \qquad k=1,2,\ldots,n-1.
\end{multline*}
The second integral is zero for $k=1,2,\ldots,n$ because of the
orthogonality condition (\ref{eq:type2ortho2}), hence we conclude that
\[ \int_0^\infty [p_{n,n}^\alpha(x)]' x^{k+\alpha} \rho_{\nu+1}(x)\, dx
   = 0, \qquad k=1,2,\ldots,n-1, \]
or equivalently
\begin{equation}  \label{eq:der1}
  \int_0^\infty [p_{n,n}^\alpha(x)]' x^{\ell+\alpha+1} \rho_{\nu+1}(x)\, dx
   = 0, \qquad \ell=0,1,\ldots,n-2.
\end{equation}
Next we do a similar analysis with (\ref{eq:type2ortho2}), in which
we use (\ref{eq:rhoup}) to find
\[  -\int_0^\infty p_{n,n}^\alpha(x) x^{k+\alpha+\nu+1} [x^{-\nu} 
 \rho_\nu(x)]'\, dx = 0, \qquad k=0,1,\ldots,n-1. \]
Integration by parts gives
\begin{multline*}
 \left. -p_{n,n}^\alpha(x) x^{k+\alpha+1} \rho_\nu(x) \right|_0^\infty \\
   +\ \int_0^\infty \left[ p_{n,n}^\alpha(x) x^{k+\alpha+\nu+1} \right]'
     x^{-\nu}\rho_\nu(x)\, dx = 0, \qquad k=0,1,\ldots,n-1. 
\end{multline*}
The integrated terms vanish for $k\geq 0$ and $\alpha > -1$. Working out
the derivative of the product then gives
\begin{multline*}
  \int_0^\infty [p_{n,n}^\alpha(x)]' x^{k+\alpha+1} \rho_\nu(x)\, dx \\
   +\ (k+\alpha+\nu+1) \int_0^\infty p_{n,n}^\alpha(x) x^{k+\alpha}
   \rho_\nu(x)\, dx = 0, \qquad k=0,1,\dots,n-1.
\end{multline*}
The second integral is zero for $k=0,1,\ldots,n-1$ because of the
orthogonality (\ref{eq:type2ortho1}), hence we have
\begin{equation}  \label{eq:der2}
  \int_0^\infty [p_{n,n}^\alpha(x)]' x^{k+\alpha+1} \rho_\nu(x)\, dx = 0,
  \qquad k=0,1,\ldots,n-1. 
\end{equation}
Now $[p_{n,n}^\alpha]'$ is a polynomial of degree $2n-1$ with leading
coefficient $2n$ (since we normalized the type 2 multiple orthogonal
polynomials by taking monic polynomials), and by (\ref{eq:der1}) and
(\ref{eq:der2}) it satisfies the orthogonality conditions
(\ref{eq:type2ortho1})--(\ref{eq:type2ortho2}) for the type 2 multiple
orthogonal polynomial $p_{n,n-1}^{\alpha+1}$. By unicity we therefore have
$[p_{n,n}^\alpha(x)]' = 2np_{n,n-1}^{\alpha+1}(x)$. 
A similar reasoning also gives the result $[p_{n,n-1}^\alpha(x)]' =
(2n-1)p_{n-1,n-1}^{\alpha+1}(x)$. Note however that the analysis does not
work when $m \notin \{n,n-1\}$.
\qed

There is a similar differential property for type 1 multiple orthogonal polynomials which complements the differential property of the type
2 multiple orthogonal polynomials given in the previous theorem.

\begin{theorem}  \label{thm:type1}
For the type 1 multiple orthogonal polynomials we have for every
$\alpha > 0$
\begin{equation}  \label{eq:type1dif}
 \frac{d}{dx} [x^\alpha q_{n,n}^\alpha(x)] = 
  x^{\alpha-1} q_{n+1,n}^{\alpha-1}(x), \quad
 \frac{d}{dx} [x^\alpha q_{n,n-1}^\alpha(x)] = 
  x^{\alpha-1} q_{n,n}^{\alpha-1}(x). 
\end{equation}
\end{theorem}
This means that the differential operator acting on $x^\alpha$ times
the type 1 polynomials raises the degree by one and lowers the $\alpha$
by one. This should be compared with the corresponding differential
property
\[   \frac{d}{dx} [ e^{-x}x^\alpha L_n^\alpha(x)]
     = (n+1) e^{-x}x^{\alpha-1} L_{n+1}^{\alpha-1}(x) \]
for Laguerre polynomials. What we will really prove is that the
derivative of $x^{\alpha} q_{n,n}^\alpha(x)$ is proportional
to $x^{\alpha-1} q_{n+1,n}^{\alpha-1}(x)$ (and similarly for the derivative
of $x^\alpha q_{n,n-1}^\alpha(x)$ which is proportional to
$x^{\alpha-1} q_{n,n}^{\alpha-1}(x)$). The choice of the proportionality
factor one in (\ref{eq:type1dif}) will fix the normalization which was
left unspecified by the homogeneous system of equations (\ref{eq:type1ortho}).

{\bf Proof:}
We begin by using (\ref{eq:type1ortho}), which we can write as
\[  \int_0^\infty x^\alpha q_{n,n}^\alpha (x^{k+1})' \, dx = 0, \qquad
  k=0,1,\ldots,2n.  \]
Integration by parts then gives
\[  \left. q_{n,n}^\alpha(x) x^{k+\alpha+1} \right|_0^\infty
    - \int_0^\infty [x^\alpha q_{n,n}^\alpha(x)]' x^{k+1}\, dx = 0, \qquad
   k=0,1,\ldots,2n. \]
The integrated terms will vanish for every $\alpha > -1$ and $k \geq 0$.
For the integrand of the integral we have
\begin{eqnarray*}
 [x^\alpha q_{n,n}^\alpha(x)]' & =  & [x^\alpha       A_{n,n}^\alpha(x)\rho_\nu(x)]'
   + [x^\alpha B_{n,n}^\alpha(x)\rho_{\nu+1}(x)]' \\
  & = & \alpha x^{\alpha-1} [A_{n,n}^\alpha \rho_\nu(x) +  B_{n,n}^\alpha(x) \rho_{\nu+1}(x)] \\
 & & +\ x^\alpha [ (A_{n,n}^\alpha(x))'\rho_\nu(x) + A_{n,n}^\alpha(x) \rho_\nu'(x) \\
 & & +\ (B_{n,n}^\alpha(x))'\rho_{\nu+1}(x) + B_{n,n}(x)\rho_{\nu+1}'(x)]. 
\end{eqnarray*}
Now use (\ref{eq:rhoup})--(\ref{eq:rhodown}) to find
\begin{eqnarray*}
   [x^\alpha q_{n,n}^\alpha(x)]' & = &
   x^{\alpha-1} [ ((\alpha+\nu) A_{n,n}^\alpha(x) + x (A_{n,n}^\alpha(x))' - xB_{n,n}^\alpha(x)) \rho_\nu(x) \\
 & &  +\ (\alpha B_{n,n}^\alpha(x) -
 A_{n,n}^\alpha(x) + x(B_{n,n}^\alpha(x))')\rho_{\nu+1}(x)] .  
\end{eqnarray*}
Observe that $(\alpha+\nu) A_{n,n}^\alpha(x) + x (A_{n,n}^\alpha(x))' - xB_{n,n}^\alpha(x)$ is a polynomial of degree at most $n+1$ and
$\alpha B_{n,n}^\alpha(x) - A_{n,n}^\alpha(x) + x(B_{n,n}^\alpha(x))'$
is a polynomial of degree at most $n$, so that we can write
\begin{equation}  \label{eq:qder1}
  [x^\alpha q_{n,n}^\alpha(x)]' = x^{\alpha-1} [P_{n+1}(x) \rho_\nu(x)
    + Q_n(x) \rho_{\nu+1}(x)] 
\end{equation}
with polynomials $P_{n+1}$ and $Q_n$ of degree at most $n+1$ and $n$
respectively. We therefore have
\[  \int_0^\infty x^{\alpha-1} [P_{n+1}(x) \rho_\nu(x)
    + Q_n(x) \rho_{\nu+1}(x)] x^{k+1}\, dx = 0, \qquad k=0,1,\ldots,2n. \]
In addition to this, we also have
\[ \int_0^\infty [x^{\alpha} q_{n,n}^\alpha(x)]' \, dx = \left. x^\alpha q_{n,n}^\alpha(x)\right|_0^\infty = 0, \]
so that
\[  \int_0^\infty [x^{\alpha} q_{n,n}^\alpha(x)]' x^{\ell}\, dx 
= 0, \qquad \ell=0,1,\ldots,2n+1. \]
In view of (\ref{eq:type1ortho}) and (\ref{eq:qder1}) this means that
\[  [x^\alpha q_{n,n}^\alpha(x)]' = \mathrm{constant\ }
   x^{\alpha-1} q_{n+1,1}^{\alpha-1}(x). \]
Since $q_{n,n}^\alpha$ and $q_{n+1,n}^{\alpha-1}$ are only
determined up to a normalizing factor, we can choose the constant
equal to one, thereby fixing the normalization.

We can repeat the analysis for $q_{n,n-1}^\alpha$ with minor changes.
In this case we have
\begin{eqnarray*}
   [x^\alpha q_{n,n-1}^\alpha(x)]' & = &
   x^{\alpha-1} [ ((\alpha+\nu) A_{n,n-1}^\alpha(x) + x (A_{n,n-1}^\alpha(x))' - xB_{n,n-1}^\alpha(x)) \rho_\nu(x) \\
 & &  +\ (\alpha B_{n,n-1}^\alpha(x) -
 A_{n,n-1}^\alpha(x) + x(B_{n,n-1}^\alpha(x))')\rho_{\nu+1}(x)] .  
\end{eqnarray*}
and since $(\alpha+\nu) A_{n,n-1}^\alpha(x) + x (A_{n,n-1}^\alpha(x))' - xB_{n,n-1}^\alpha(x)$ is a polynomial of degree at most $n$ and
$\alpha B_{n,n-1}^\alpha(x)-A_{n,n-1}^\alpha(x) + x(B_{n,n-1}^\alpha(x))'$
a polynomial of degree at most $n$ also, the orthogonality
\[  \int_0^\infty [x^\alpha q_{n,n-1}^\alpha(x)]' x^{\ell} \, dx = 0,
 \qquad \ell = 0,1,\ldots,2n \]
gives the required result.

Observe that the reasoning does not work for $[x^\alpha q_{n,m}^\alpha(x)]'$
when $m \notin \{n,n-1\}$.
\qed

The proof also shows that
\begin{eqnarray*}
   A_{n+1,n}^{\alpha-1}(x) & = & (\alpha+\nu) A_{n,n}^\alpha(x) + 
   x [A_{n,n}^\alpha(x)]' - x B_{n,n}^\alpha(x), \\
   B_{n+1,n}^{\alpha-1}(x) & = & \alpha B_{n,n}^\alpha(x) -       A_{n,n}^\alpha(x) + x[B_{n,n}^\alpha(x)]',
\end{eqnarray*}
and similarly
\begin{eqnarray*}
   A_{n,n}^{\alpha-1}(x) & = & (\alpha+\nu) A_{n,n-1}^\alpha(x) + 
   x [A_{n,n-1}^\alpha(x)]' - x B_{n,n-1}^\alpha(x), \\
   B_{n,n}^{\alpha-1}(x) & = & \alpha B_{n,n-1}^\alpha(x) -      
   A_{n,n-1}^\alpha(x) + x[B_{n,n-1}^\alpha(x)]'.
\end{eqnarray*}  

\setcounter{equation}{0}
\section{Rodrigues formula}
As a consequence of Theorem \ref{thm:type1} we have the following Rodrigues formula for the type 1 multiple orthogonal polynomials.

\begin{theorem}  \label{thm:rodrigues}
The type 1 multiple orthogonal polynomials can be obtained from 
\begin{equation}  \label{eq:rodrigues}
  \frac{d^{2n}}{dx^{2n}} \left( x^{2n+\alpha} \rho_\nu(x) \right) =
     x^\alpha q_{n,n-1}^\alpha(x), \quad
  \frac{d^{2n+1}}{dx^{2n+1}} \left( x^{2n+1+\alpha} \rho_\nu(x) \right) =
     x^\alpha q_{n,n}^\alpha(x).
\end{equation}
\end{theorem}

{\bf Proof:}
By combining the two formulas in (\ref{eq:type1dif}) we get
\[  \frac{d^2}{dx^2} \left( x^\alpha q_{n,n-1}^\alpha(x) \right)
  = x^{\alpha-2} q_{n+1,n}^{\alpha-2}. \]
Iterate this $k$ times to get
\[   \frac{d^{2k}}{dx^{2k}} \left( x^\alpha q_{n,n-1}^\alpha(x) \right)
  = x^{\alpha-2k} q_{n+k,n+k-1}^{\alpha-2k}. \]
Choose $n=0$ and $\alpha=2k+\beta$ to find
\[    \frac{d^{2k}}{dx^{2k}} \left( x^{2k+\beta} q_{0,-1}^{2k+\beta}(x) \right)
  = x^{\beta} q_{k,k-1}^{\beta}. \]
Now $q_{0,-1}^\alpha(x) = A_{0,-1}^\alpha \rho_\nu(x)$, where $A_{0,-1}^\alpha$ is a constant, which we will take equal to one to fix
normalization. Thus, if we replace $\beta$ by $\alpha$ and $k$ by $n$,
then we find the first formula in (\ref{eq:rodrigues}).

Similarly, we can use
\[   \frac{d}{dx} \left( x^\alpha q_{n,n-1}^\alpha(x) \right)
    = x^{\alpha-1} q_{n,n}^{\alpha-1}(x), \]
and differentiate it $2k$ times, to find
\[  \frac{d^{2k+1}}{dx^{2k+1}} \left( x^\alpha q_{n,n-1}^\alpha(x) \right)
  = x^{\alpha-2k-1} q_{n+k,n+k}^{\alpha-2k-1}(x). \]
Choose $n=0$ and $\alpha=2k+1+\beta$ to find
\[  \frac{d^{2k+1}}{dx^{2k+1}} \left( x^{2k+1+\beta} \rho_\nu(x)\right)
  = x^{\beta} q_{k,k}^{\beta}(x). \]
Replacing $\beta$ by $\alpha$ and $k$ by $n$ then gives the required formula. \qed

The Rodrigues formula allows us to compute the type 1 multiple
orthogonal polynomials explicitly. There is a relationship
between type 1 and type 2 polynomials (which is not typical for the Macdonald weights but holds in general) and this allows us to compute
the type 2 multiple polynomials as well. Indeed, we can write
\[  \begin{pmatrix}
       q_{n,n}^\alpha \\  q_{n,n-1}^\alpha 
    \end{pmatrix} =
    \begin{pmatrix}
     A_{n,n}^\alpha  & B_{n,n}^\alpha  \\
     A_{n,n-1}^\alpha  & B_{n,n-1}^\alpha 
    \end{pmatrix}
     \begin{pmatrix}
      \rho_\nu  \\  \rho_{\nu+1} 
     \end{pmatrix} ,  \]
so that we have
\[  \begin{pmatrix}
     A_{n,n}^\alpha  & B_{n,n}^\alpha  \\
     A_{n,n-1}^\alpha  & B_{n,n-1}^\alpha 
    \end{pmatrix}^{-1}
    \begin{pmatrix}
       q_{n,n}^\alpha  \\  q_{n,n-1}^\alpha 
    \end{pmatrix} =
    \begin{pmatrix}
      \rho_\nu  \\  \rho_{\nu+1} 
     \end{pmatrix}.  \]
Writing the inverse of a matrix as the adjoint matrix divided by the
determinant gives
\begin{equation}  \label{eq:matrix}
  \begin{pmatrix}
     B_{n,n-1}^\alpha  & -B_{n,n}^\alpha  \\
     -A_{n,n-1}^\alpha  & A_{n,n}^\alpha 
    \end{pmatrix}
    \begin{pmatrix}
       q_{n,n}^\alpha  \\  q_{n,n-1}^\alpha 
    \end{pmatrix} =
     [A_{n,n}^\alpha B_{n,n-1}^\alpha  - 
    A_{n,n-1}^\alpha B_{n,n}^\alpha ]
    \begin{pmatrix}
      \rho_\nu  \\  \rho_{\nu+1} 
     \end{pmatrix}.  
\end{equation}
The polynomial $A_{n,n}^\alpha B_{n,n-1}^\alpha  - 
    A_{n,n-1}^\alpha B_{n,n}^\alpha $ is of degree at most $2n$
and satisfies
\begin{multline*}
 \int_0^\infty [A_{n,n}^\alpha(x) B_{n,n-1}^\alpha(x) - 
    A_{n,n-1}^\alpha(x)B_{n,n}^\alpha(x)] \rho_\nu(x) x^{k+\alpha} \, dx \\
   = \int_0^\infty [B_{n,n-1}^\alpha(x) q_{n,n}^\alpha(x)
   - B_{n,n}^\alpha(x) q_{n,n-1}^\alpha(x)] x^{k+\alpha} \, dx
  = 0, \qquad k=0,1,\ldots,n-1, 
\end{multline*}
where we have used (\ref{eq:matrix}) and (\ref{eq:type1ortho}).
Furthermore we also have
\begin{multline*}
 \int_0^\infty [A_{n,n}^\alpha(x)B_{n,n-1}^\alpha(x) - 
    A_{n,n-1}^\alpha(x)B_{n,n}^\alpha(x)] \rho_{\nu+1}(x) x^{k+\alpha} \, dx \\
   = \int_0^\infty [-A_{n,n-1}^\alpha(x) q_{n,n}^\alpha(x)
   + A_{n,n}^\alpha(x) q_{n,n-1}^\alpha(x)] x^{k+\alpha} \, dx
  = 0, \qquad k=0,1,\ldots,n-1. 
\end{multline*}
Hence we can conclude that
\[   A_{n,n}^\alpha(x)B_{n,n-1}^\alpha(x) - 
    A_{n,n-1}^\alpha(x)B_{n,n}^\alpha(x) = \mathrm{constant\ } p_{n,n}^\alpha(x), \]
where the constant is the leading coefficient of the polynomial.
A similar reasoning using
\[ \begin{pmatrix}
       q_{n+1,n}^\alpha  \\  q_{n,n}^\alpha 
    \end{pmatrix} =
    \begin{pmatrix}
     A_{n+1,n}^\alpha  & B_{n+1,n}^\alpha  \\
     A_{n,n}^\alpha  & B_{n,n}^\alpha 
    \end{pmatrix}
     \begin{pmatrix}
      \rho_\nu  \\  \rho_{\nu+1} 
     \end{pmatrix} ,  \]
gives
\[  A_{n+1,n}^\alpha(x) B_{n,n}^\alpha(x) - A_{n,n}^\alpha(x)     B_{n+1,n}^\alpha(x) = \mathrm{constant\ } p_{n+1,n}^\alpha(x). \]

\setcounter{equation}{0}
\section{Recurrence relation}
To simplify the notation we put
\[  P_{2n}(x) = p_{n,n}^\alpha(x), \quad P_{2n+1}(x) = p_{n+1,n}^\alpha(x). \]
It is known that the sequence $P_n(x)$, $n=0,1,2,\ldots$ satisfies a
third order recurrence relation of the form
\begin{equation}  \label{eq:recur}
  xP_n(x) = P_{n+1}(x) + b_n P_n(x) + c_n P_{n-1}(x) + d_n P_{n-2}(x),
\end{equation} 
(see, e.g., \cite{kalia} or \cite{nikisor}). Explicit formulas for the recurrence coefficients are given in the following theorem.

\begin{theorem}  \label{thm:recurrence}
The recurrence coefficients in (\ref{eq:recur}) are given by
\begin{eqnarray*}
   b_n & = & (n+\alpha+1)(3n+\alpha+2\nu) - (\alpha+1)(\nu-1) \\
   c_n & = & n(n+\alpha)(n+\alpha+\nu)(3n+2\alpha+\nu) \\
   d_n & = & n(n-1)(n+\alpha-1)(n+\alpha)(n+\alpha+\nu-1)(n+\alpha+\nu).
\end{eqnarray*}
\end{theorem}

{\bf Proof:}
We begin with the recurrence coefficients of even index,
which are used in the recurrence relation 
\begin{equation}  \label{eq:even}
 xp_{n,n}^\alpha(x) = p_{n+1,n}^\alpha(x) + b_{2n}p_{n,n}^\alpha(x) + c_{2n} p_{n,n-1}^\alpha(a)
 + d_{2n} p_{n-1,n-1}^\alpha(x). 
\end{equation}
We will write the polynomials $p_{n,m}^\alpha$ explicitly
as
\[   p_{n,m}^\alpha(x) = \sum_{k=0}^{n+m} a_{n,m}^\alpha(k) x^{n+m-k}, \]
and since we are dealing with monic polynomials, we have
\[   a_{n,m}^\alpha(0) = 1. \]
Comparing the coefficient of $x^{2n}$ in (\ref{eq:even}) gives
\begin{equation}  \label{eq:b2n}
  b_{2n} = a_{n,n}^\alpha(1) - a_{n+1,n}^\alpha(1),
\end{equation}
hence we need to know the coefficients $a_{n,m}^\alpha(1)$. Comparing
the coefficient of $x^{2n-2}$ and $x^{2n-3}$ respectively in (\ref{eq:type2dif}) gives the recurrence
\[   (2n-1)a_{n,n}^\alpha(1) = 2n a_{n,n-1}^{\alpha+1}(1), \quad
     (2n-2) a_{n,n-1}^\alpha(1) = (2n-1) a_{n-1,n-1}^{\alpha+1}(1). \]
Combining these  relations gives
\[     a_{n,n}^\alpha(1) = \frac{n}{n-1} a_{n-1,n-1}^{\alpha+2}(1), \]
which leads to
\[  a_{n,n}^\alpha(1) = n a_{1,1}^{\alpha+2n-2}(1). \]
In order to obtain an explicit formula, we compute $p_{1,1}^\alpha$
explicitly by solving the system of equations
\begin{eqnarray*}
   \int_0^\infty (x^2 + a_{1,1}^\alpha(1) x + a_{1,1}^\alpha(2)) x^\alpha
  \rho_\nu(x) \, dx & = &  0 \\
   \int_0^\infty (x^2 + a_{1,1}^\alpha(1) x + a_{1,1}^\alpha(2)) x^\alpha
  \rho_{\nu+1}(x) \, dx & = &  0.
\end{eqnarray*}
Using the moments (\ref{eq:rhomoment}) this gives
\begin{eqnarray}
   a_{1,1}^\alpha(1) & = &  -2(2+\alpha)(2+\alpha+\nu), \label{eq:a111} \\
   a_{1,1}^\alpha(2) & = & (1+\alpha)(2+\alpha)(1+\alpha+\nu)(2+\alpha+\nu).  \label{eq:a112}
\end{eqnarray}		
This gives
\begin{equation}  \label{eq:ann1}
     a_{n,n}^\alpha(1) = -2n(\alpha+2n)(\alpha+2n+\nu),
\end{equation}
and since $2na_{n+1,n}^\alpha(1) = (2n+1) a_{n,n}^{\alpha+1}(1)$ this also gives
\begin{equation}  \label{eq:an+1n1}
   a_{n+1,n}^\alpha(1) = -(2n+1)(\alpha+2n+1)(\alpha+2n+\nu+1).
\end{equation}
Inserting this in (\ref{eq:b2n}) gives the requested formula for the recurrence coefficient $b_{2n}$.

Next we compare the coefficient of $x^{2n-1}$ in (\ref{eq:even}) to find
\begin{equation}  \label{eq:c2n}
  c_{2n} = a_{n,n}^\alpha(2) - a_{n+1,n}^\alpha(2) - b_{2n} a_{n,n}^\alpha(1). 
\end{equation}
This means that we also need to know $a_{n,n}^\alpha(2)$ and 
$a_{n+1,n}^\alpha(2)$.  Compare
the coefficient of $x^{2n-3}$ and $x^{2n-4}$ respectively in (\ref{eq:type2dif}), then
\[  (2n-2)a_{n,n}^\alpha(2) = 2n a_{n,n-1}^{\alpha+1}(2), \quad
    (2n-3)a_{n,n-1}^\alpha(2) = (2n-1) a_{n-1,n-1}^{\alpha+1}(2), \]
which combined gives
\[   a_{n,n}^\alpha(2) = \frac{(2n)(2n-1)}{(2n-2)(2n-3)}
    a_{n-1,n-1}^{\alpha+2}(2) = \frac{(2n)(2n-1)}{2}
  a_{1,1}^{\alpha+2n-2}(2). \]
Using (\ref{eq:a112}) this gives
\begin{equation}  \label{eq:ann2}
  a_{n,n}^\alpha(2) = n(2n-1)(\alpha+2n-1)(\alpha+2n)(\alpha+2n+\nu-1)
   (\alpha+2n+\nu),
\end{equation}
and since $(2n-1)a_{n+1,n}^\alpha(2) = (2n+1) a_{n,n}^{\alpha+1}(2)$
this also gives
\begin{equation}  \label{eq:an+1n2}
  a_{n+1,n}^\alpha(2) = n(2n+1)(\alpha+2n)(\alpha+2n+1)(\alpha+2n+\nu)
   (\alpha+2n+\nu+1).
\end{equation}
Inserting (\ref{eq:ann2}), (\ref{eq:an+1n2}), (\ref{eq:ann1}) and
the expression for $b_{2n}$ into (\ref{eq:c2n}) gives, after some
straightforward calculus (or by using Maple) the requested expression for
$c_{2n}$.

Finally, compare the coefficient of $x^{2n-2}$ in (\ref{eq:even}), then
\begin{equation}  \label{eq:d2n}
  d_{2n} = a_{n,n}^\alpha(3) - a_{n+1,n}^\alpha(3) - 
  b_{2n} a_{n,n}^\alpha(2) - c_{2n} a_{n,n-1}^\alpha(1), 
\end{equation}
so that we need $a_{n,n}^\alpha(3)$ and $a_{n+1,n}^\alpha(3)$.
To this end we compare the coefficient of $x^{2n-4}$ and $x^{2n-5}$
respectively in (\ref{eq:type2dif}) to find
\[ (2n-3) a_{n,n}^\alpha(3) = 2n a_{n,n-1}^\alpha(3), \quad
   (2n-4) a_{n,n-1}^\alpha(3) = (2n-1) a_{n-1,n-1}^\alpha(3), \]
which combined gives
\[  a_{n+1,n}^\alpha(3) = \frac{2n(2n+1)}{(2n-2)(2n-3)} 
   a_{n,n-1}^{\alpha+2}(3) = \frac{(2n+1)(2n)(2n-1)}{6} 
 a_{2,1}^{\alpha+2n-2}(3). \]
In order to find $a_{2,1}^\alpha(3)$ we will compute $p_{2,1}^\alpha$
explicitly by solving the system of equations
\begin{eqnarray*}
 \int_0^\infty (x^3 + a_{2,1}^\alpha(1) x^2 + a_{2,1}^\alpha(2) x
   + a_{2,1}^\alpha(3)) x^\alpha \rho_\nu(x) \, dx & = & 0 \\
 \int_0^\infty (x^3 + a_{2,1}^\alpha(1) x^2 + a_{2,1}^\alpha(2) x
   + a_{2,1}^\alpha(3)) x^{\alpha+1} \rho_\nu(x) \, dx & = & 0 \\
  \int_0^\infty (x^3 + a_{2,1}^\alpha(1) x^2 + a_{2,1}^\alpha(2) x
   + a_{2,1}^\alpha(3)) x^\alpha \rho_{\nu+1}(x) \, dx & = & 0
\end{eqnarray*}
which by using (\ref{eq:rhomoment}) and some calculus gives
\begin{equation}  \label{eq:a213}
  a_{2,1}^\alpha(3) = -(3+\alpha+\nu)(2+\alpha+\nu)(1+\alpha+\nu)
  (3+\alpha)(2+\alpha)(1+\alpha).
\end{equation}
Using this gives
\begin{multline} \label{eq:an+1n3}
  a_{n+1,n}^\alpha(3) = - \frac{(2n+1)(2n)(2n-1)}{6} \\
   (\alpha+2n+\nu+1)(\alpha+2n+\nu)(\alpha+2n+\nu-1)
   (\alpha+2n+1)(\alpha+2n)(\alpha+2n-1), 
\end{multline}
and since $(2n-3)a_{n,n}^\alpha(3) = 2n a_{n,n-1}^{\alpha+1}(3)$ we also have
\begin{multline}  \label{eq:ann3}
  a_{n,n}^\alpha(3) = - \frac{2n(2n-1)(2n-2)}{6}   \\
   (\alpha+2n+\nu)(\alpha+2n+\nu-1)(\alpha+2n+\nu-2)
   (\alpha+2n)(\alpha+2n-1)(\alpha+2n-2).
\end{multline}
Using (\ref{eq:ann3}), (\ref{eq:an+1n3}), (\ref{eq:ann2}), (\ref{eq:an+1n1})
and the formulas for $b_{2n}$ and $c_{2n}$ in (\ref{eq:d2n})
then gives the requested expression for $d_{2n}$.

In a similar way we proceed with the odd indices which appear in the recurrence relation
\begin{equation} \label{eq:odd}
  xp_{n+1,n}^\alpha(x) = p_{n+1,n+1}^\alpha(x) + 
   b_{2n+1} p_{n+1,n}^\alpha(x) + c_{2n+1} p_{n,n}^\alpha(x)
  + d_{2n+1} p_{n,n-1}^\alpha(x).
\end{equation}
By comparing the coefficients of $x^{2n+1}$, $x^{2n}$ and $x^{2n-1}$
we get expressions for $b_{2n+1}$, $c_{2n+1}$ and $d_{2n+1}$
respectively in terms of the coefficients $a_{n,m}^\alpha$ and after
working out these expressions we get the required formulas. \qed

Observe that the recurrence coefficients have the asymptotic behavior
\begin{equation}  \label{eq:recasym}
   b_n \sim 3n^2, \quad c_n \sim 3 n^4, \quad c_n \sim n^6. 
\end{equation}

\bigskip

\begin{verbatim}
Department of Mathematics
Katholieke Universiteit Leuven
Celestijnenlaan 200 B
B-3001 Heverlee (Leuven)
Belgium
\end{verbatim}

\end{document}